\documentclass{amsart}
\usepackage{amsmath,amsthm,amssymb,amsfonts}

\usepackage{graphicx}
\usepackage{tikz}

\usepackage{epsfig}
\usepackage[pagebackref, colorlinks=true]{hyperref}

\hypersetup{urlcolor=blue, citecolor=blue}

\def\doi#1{   {\href{http://dx.doi.org/#1}
		{{\mdseries\ttfamily DOI}}}}

\usepackage{graphics}

\def\p{\partial}

\newcommand{\al}{\alpha}    \newcommand{\be}{\beta}
\newcommand{\de}{\delta}    
  \newcommand{\ep}{\varepsilon}

\newcommand{\la}{\lambda}

\newcommand{\R}{\mathbb{R}}\newcommand{\Z}{\mathbb{Z}}
\newcommand{\N}{\mathbb{N}}\newcommand{\C}{\mathbb{C}}

\newcommand{\half}{\frac{1}{2}}

\newcommand{\pt}{\partial_t}

\newcommand{\les}{{\lesssim}}
\newcommand{\beeq}{\begin{equation}}\newcommand{\eneq}{\end{equation}}

\newcommand{\supp}{\text{supp}}

\newenvironment{prf}{\noindent {\bf Proof.} }{\endprf\par}
\def \endprf{\hfill  {\vrule height6pt width6pt depth0pt}\medskip}

\def\<{\langle}             \def\>{\rangle}
\def\({\left(}                 \def\){\right)}
\numberwithin{equation}{section}
\newtheorem{thm}{Theorem}[section]
\newtheorem{cor}[thm]{Corollary}
\newtheorem{lem}[thm]{Lemma}
\newtheorem{prop}[thm]{Proposition}
\newtheorem{defn}[thm]{Definition}
\newtheorem{rem}[thm]{Remark}

\title[semilinear wave equations with 
energy supercritical powers]
{Global solutions with small initial data to semilinear wave equations with energy supercritical powers}

\author{Kerun Shao}
\address{School of Mathematical Sciences\\ Zhejiang University\\Hangzhou 310058, P. R. China}\email{shaokr@163.com}

\author{Chengbo Wang$^{*}$}\thanks{* Corresponding author}
\address{School of Mathematical Sciences\\ Zhejiang University\\Hangzhou 310058, P. R. China}\email{wangcbo@zju.edu.cn}



\date{\today}

\thanks{
	The authors were supported  by NSFC 
	11971428 and  NSFC 12141102.}

\begin{document}
	\bibliographystyle{plain}
	\maketitle

	\begin{abstract}
		Considering $1+n$ dimensional semilinear wave equations with 
		energy supercritical powers $p> 1+4/(n-2)$, 
		we obtain global solutions for any 
		initial data with small norm in  $H^{s_c}\times H^{s_c-1}$,
		under the technical smooth condition $p>s_c-\bar{s}_0$, 
		with $\bar{s}_0= 1/2+(n-3)/(2\max(n-1-p,n-3))$ and $s_c=n/2-2/(p-1)$.
		In particular, combined with previous works, our results give a complete verification of the Strauss conjecture, up to space dimension $9$. The higher dimensional case, $n\ge 10$, seems to be unreachable, in view of the wellposed theory in $H^s$.
	\end{abstract}
	
\section{Introduction}
Let $n\ge 2$, $F_p(u)=|u|^p$ with 
$p>1$.
Given $u_0, u_1 \in C_c^{\infty}(\R^n)$, 
considering the following Cauchy problem with sufficiently small initial data (of size $\ep >0$):
\begin{equation}\label{eq-Strauss}
	\left\{
	\begin{aligned}
		&\Box u:=(\pt^2-\Delta)u=
		F_p(u)\ ,\\
		&u(0)=\ep u_0\ ,\ \partial_t u(0)=\ep u_1\ ,
	\end{aligned}
	\right.
\end{equation}
the Strauss conjecture claims that the critical power, denoted by $p_c(n)$, for the global existence vs. blow-up is given by the positive root of the quadratic equation
$(n-1)p^2-(n+1)p-2=0$.
The problem has been well investigated and the 
conjecture is well believed to be verified for all spatial dimensions.
However, if we take a closer look at the results obtained in literature, we find a slight gap. More precisely, the conjecture has been verified for any sub-conformal powers, $1<p\le p_{conf}:=1+4/(n-1)$ (noticing that $p_c\in (1, p_{conf})$), see \cite{John79}, \cite{Glassey81bu},
\cite{Sideris84} for the blow-up results with $p\in (1, p_c)$, and
\cite{John79}, \cite{Glassey81ex},
\cite{Zhou95}, \cite{LdSo96}, \cite{GLS97}, \cite{Ta02} for global results with $p\in (p_c, p_{conf}]$.
See also \cite{Wang18} for recent review.

Turning to the super-conformal powers, let us give a quick review on the closely related well-posed theory. The problem \eqref{eq-Strauss} is scale-invariant, 
a solution $u$ will give rise to a class of solutions $$u_\la(t,x)=\la^{-2/(p-1)}u(t/\la, x/\la) \ , \  \la>0\ ,$$
which gives a lower bound of the regularity $s$ so that the problem is locally well-posed in $C_t H^s\cap C^1_t H^{s-1}$, that is,
\beeq\label{eq-sc}
s\ge s_c:=\frac n 2-\frac{2}{p-1}\ ,\eneq
see, e.g., \cite{LdSo95}.
Heuristically, 
if we look for
$C_t {H} ^{s_c}\cap
C_t^1 H ^{s_c-1}$ solutions,
in view of the energy estimates,
a necessary least regularity requirement for the nonlinear term would be
$p\ge s_c-1$,
so that as a function of the unknown function, $F(u)=|u|^p\in C^{s_c-1}$.
In this sense, 
in purpose of the wellposedness in $H^{s_c}$,
it is then natural to pose the following restriction (up to the critical case) on the smoothness of the nonlinearity
\beeq \label{eq-regu}p> s_c-1\ .\eneq

Observing that such a restriction, \eqref{eq-regu}, takes effect only when $n\ge 10$, in which case,
there always exists  $p>2$ such that $s_c-1>p$. For example, when  $n=10$,
we have
$$s_c-1>p\Leftrightarrow
5-2/(p-1)-1>p
\Leftrightarrow
(p-1)^2-3(p-1)+2<0
\Leftrightarrow p\in (2,3)
\ .
$$
Then, it is natural to infer that the Strauss conjecture may hold for spatial dimension up to $n=9$, while for the high dimensional situation $n\ge 10$, we may need to develop alternative methods beyond the well-posed theory in $C_t H^s$.

Now, let us turn to the current state of art for the Strauss conjecture with super-conformal powers, $p\ge p_{conf}$.
At first, when $n\le 4$, it is well-known that the Strauss conjecture is valid for any $p>p_c$, partly due to the fact that $p_c\ge 2$, see, e.g., \cite{FW11}, \cite{HMSSZ}. 
For $n\ge 5$, the well-posed theory developed in
\cite{LdSo95}, \cite{NO99} is sufficient to give critical well-posed result, which, in particular concludes global existence under the technical requirement that
$F_p$ is algebraic, or
$$p>[s_c]\ ,$$
where $[s_c]$ denotes the integer part of $s_c$.
In particular, when $s_c\le 1$,
i.e., for energy subcritical or critical powers
$$p\le p_{H^1}:=1+\frac{4}{n-2}\ ,$$
we have global existence with small data.
Through a direct calculation, we see that we have had a complete proof of the Strauss conjecture up to spatial dimension $n=7$.

When $n=8$,
we have $p>[s_c]$, for all $p>1$ except the case $p=3$.
For the excluded case with $p=3$ and $F_p(u)=|u|^3$, we have $s_c=3=p$ and the requirement $p>[s_c]$ barely hold. 
The serious obstacle appears when $n\ge 9$. For example,
when $n=9$ and $7/3<p< 3$, we have $s_c= n/2-2/(p-1)\in (3, 7/2)$, and so $[s_c]=3>p$.

The aim of this paper is to shed new light on this classical problem, by trying to prove global existence for more powers. In particular, we would like to see if the Strauss conjecture holds for $n=9$ or not.

As we have recalled, the Strauss conjecture has been verified for any energy subcritical and critical powers. In view of this fact, 
here and in what follows,
we will restrict ourselves to the case of energy supercritical powers, $p>p_{H^1}$.
Now, we are ready to present our main result,
see Theorem \ref{thm-mainv2} for a more precise statement once further notations are established.

\begin{thm}\label{thm-main}
	Let $4\le n\le 9$, $p>p_{H^1}$.
	Then, for any $u_0, u_1\in C_c^\infty(\R^n)$,
	there exists $\ep_0>0$, such that
	the problem \eqref{eq-Strauss} admits a global weak solution, 
	as long as $\ep\in (0,
	\ep_0)$. In addition, when $n\ge 10$, the same result applies for any
	$p>p_{H^1}$ with $p>s_c-\bar{s}_0$, where \beeq\label{eq-s0}\bar{s}_0=\half+\frac{n-3}{2\max(n-1-p,n-3)}\ .\eneq 
	To be more specific, the admissible range of $p$ for $n\ge 10$ is $p\in (p_{H^1}, a)\cup (b,\infty)$, where $a$ is the smallest real root of the equation $a^3-\frac{3n-1}{2}a^2+\frac{n^2+5}{2}a-\frac{n^2+n}{2}=0$ and 
	$b=\frac{n+\sqrt{(n-4)^2-32}}{4}$.
\end{thm}

In particular, we see that the Strauss conjecture holds for any spatial dimension $n\le 9$.

\begin{rem}\label{rem-main2}
	Actually, the compactness of the initial data is not necessary. The conditions on the initial data can be relaxed to the realization of homogeneous Sobolev space $(u_0, u_1)\in \dot{\mathcal{H}}^{s_c}\times \dot{\mathcal{H}}^{s_c-1}$.
	Moreover, the global solution is unique in $L_t^q(\R,L_x^r(\R^n))$, for certain
	choice of $(q,r)$ (see \eqref{eq-nlowp} and \eqref{eq-nhighp}).
	Here, we follow \cite{BCD} to define the realization of homogeneous Besov and Sobolev spaces, $\dot{\mathcal{B}}_r^s$ and $\dot{\mathcal{H}}^s$, and we recall some fundamental definitions and results in the Appendix. In addition, to obtain global solution for more powers, we give up the requirement that the solution lies in the space $C_t\dot{\mathcal{H}}^{s_c}$,
	and use some subspace of
	$C^1_t S'$ instead.
	Here and in what follows, we use $S$ and $S'$   
	to denote the space of Schwartz functions and tempered distributions.
\end{rem}
\begin{rem}\label{rem-main1}
	As is clear from the proof, under the prerequisite $s_c-\bar{s}_0<p$, the nonlinearity $F_p(u)=|u|^p$ can be replaced by any functions $f(u)$ satisfying, with $k:=[s_c-\bar{s}_0]$, \\
	\begin{enumerate}
		\item $f\in C^{k}(\C;\C)$, and $f^{(j)}(0)=0$ for all $j$ with $0\leq j \leq k$\ ,
		\item There exists a constant $C$ such that for all $z_1,z_2\in \C$
		$$|f^{k}(z_1)-f^{k}(z_2)|\leq
		\left\{
		\begin{array}{ll}
			C(|z_1|^{p-k-1}+|z_2|^{p-k-1})|z_1-z_2|\ ,& p\geq k+1\ ,\\
			C|z_1-z_2|^{p-k}\ , & p< k+1\ .
		\end{array}
		\right.
		$$
	\end{enumerate}
\end{rem}
\begin{rem}
	In \cite{NO99-2}, \cite{NO01},
	Nakamura and Ozawa have investigated the similar problem with much more general nonlinearity (which can grow exponentially near infinity), when the initial 
	data have more regularity, i.e., $(u_0, u_1)\in
	(\dot{\mathcal{H}}^{s}\cap \dot{\mathcal{H}}^{1/2})
	\times(\dot{\mathcal{H}}^{s-1}\cap \dot{\mathcal{H}}^{-1/2})$, for some $s\ge n/2$.
	It might be interesting to investigate to what extent we can adapt the method in this paper to improve the conditions.
\end{rem}

Let us conclude the introduction with  the idea of proof. In the previous works,
\cite{LdSo95}, \cite{NO99}, the authors exploited the Strichartz estimates and the fractional chain rule.
The Strichartz estimates used in these works are 
essentially those derived from the homogeneous Strichartz estimates, which were available in \cite{LdSo95} (see also \cite{GV95}, \cite{KeelTao98}).
On the one hand, it is well-known that the homogeneous Strichartz estimates could not be extended for more pairs $(q,r)$, see, e.g., \cite{KeelTao98}, \cite{FW2},
\cite{GLNY18}.
On the other hand, 
the fractional chain rule used in
\cite{NO99} seems to be sharp in general.
It is then natural to take an alternative route, by exploiting the  inhomogeneous Strichartz estimates, which are known to hold for more
general pairs $(q,r)$,
available in  
\cite{MR1046747, MR1052018, MR2134950}.
Actually, with the help of Harmse-Oberlin estimates (\cite{MR1046747,MR1052018}), one may prove the Strauss conjecture for $p>s_c-(1-1/(2n))$ and $p>p_{H^1}$, which has been sufficient to give Theorem \ref{thm-main} for $n\le 9$.

\section{Preliminary estimates}
In this section, we give the necessary linear and nonlinear estimates.

\subsection{Strichartz estimates}
Considering the linear wave equation:
$$\begin{cases}
	\Box u(t,x)=F(t,x)\ ,\\
	u(0,x)=f(x) \ , \ u_t(0,x)=g(x)\ ,
\end{cases}$$
the solution could be written as the superposition of the homogeneous part $S(t)(f, g)$ (with $F=0$)
and the inhomogeneous part $\mathcal{G}F(t)$ (with $f=g=0$):
$$u=S(t)(f, g)+\mathcal{G}F(t)\ .$$

At first, let us recall the homogeneous and inhomogeneous Strichartz estimates.
\begin{defn}\label{defn-sigama-adiacc}
	The exponent pair $(q,r)$ is $\sigma$-admissible, if $q,r\in [2,\infty]$, $(q,r,\sigma)\neq(2,\infty,1)$ and
	$$\frac{1}{q}\leq \sigma(\frac{1}{2}-\frac{1}{r})\ .$$
	The exponent pair $(q,r)$ is $\sigma$-acceptable, if $$q\in [1,\infty)\ , \ r\in [2,\infty]\ , \ \frac{1}{q}< 2\sigma(\frac{1}{2}-\frac{1}{r}) \  \mbox{\rm or } \  (q,r)=(\infty,2).$$
\end{defn}

\begin{rem}
	In the following content, we will always take $n\ge 4$ and $\sigma=\frac{n-1}{2}>1$.
\end{rem}

\begin{lem}[Homogeneous Strichartz estimates, \cite{KeelTao98}]\label{lem-hoestimate-Tao}
	If $(q,r),(\tilde{q},\tilde{r})$ are $\sigma$-admissible, then we have the estimates
	\begin{equation}
		\begin{aligned}
			&\|S(t)(f, g)\|_{L_t^q(\R;\dot{B}_r^{\bar{s}})}\lesssim \|(f,g)\|_s:=
			\|f\|_{\dot{\mathcal{H}}^{s}}+\|g\|_{\dot{\mathcal{H}}^{s-1}}
			\ ,\ s=\bar{s}+\frac{n}{2}-\frac{1}{q}-\frac{n}{r}\ ,\\
			&\|\mathcal{G}F\|_{L_t^q(\R;\dot{B}_{r}^{\bar{s}})}\lesssim \|F\|_{L_t^{\tilde{q}'}(\R;\dot{B}_{\tilde{r}'}^s)}\ ,\  s=\bar{s}+\frac{1}{\tilde{q}'}+\frac{n}{\tilde{r}'}-\frac{1}{q}-\frac{n}{r}-2\ .\nonumber
		\end{aligned}
	\end{equation}

\end{lem}

In contrast to the homogeneous estimates, the inhomogeneous estimate could be true for more general choice of pairs.
\begin{lem}[Inhomogeneous Strichartz estimates, \cite{MR2134950}]\label{lem-inhoestimate-Fos}
	Let $n\ge 4$, and $(q,r)$, $(\tilde{q},\tilde{r})$ be $\sigma$-acceptable pairs. Then  the estimate
	$$\|P_0(\mathcal{G}F)\|_{L_t^q(\R;L_{\R^n}^r)}\lesssim \|P_0 F\|_{L_t^{\tilde{q}'}(\R;L_{\R^n}^{\tilde{r}'})}$$
	holds,  provided that
	\beeq\label{eq-Fos}\frac{1}{q}+\frac{1}{\tilde{q}}=\sigma(1-\frac{1}{r}-\frac{1}{\tilde{r}})<1 \ , \
	\frac{\sigma-1}{r} \leq \frac{\sigma}{\tilde{r}}\ , \  \frac{\sigma-1}{\tilde{r}} \leq \frac{\sigma}{r} \ .
	\eneq
	Here, $P_0$ is the Littlewood-Paley projector with frequency $|\xi|\sim 1$. See the definition of $P_0$ in \cite[Section 2.2]{BCD}.
\end{lem}
\begin{rem}
	Actually, the proof of Lemma \ref{lem-inhoestimate-Fos} is established by the estimates for the operators, $\int_{0}^{t}{e^{\pm i (t-s)\sqrt{-\Delta}}\over\sqrt{-\Delta}}P_0F(s)ds$, \cite{KeelTao98}. 
	Set $$\mathcal{H}F(t)=\int_{0}^{t}{\cos (t-s)\sqrt{-\Delta}\over\sqrt{-\Delta}}F(s)ds\ .$$Lemma \ref{lem-inhoestimate-Fos}  still holds for the operator $\mathcal{H}$ instead of $\mathcal{G}$.
\end{rem}

By homogenous Littlewood-Paley decomposition and the Bernstein inequality, we get the following version of the inhomogeneous Strichartz estimates.
\begin{cor}\label{cor-inhoesti-wave}
	Let  $(q,r_\alpha)$ and $(\tilde{q},\tilde{r}_\beta)$ satisfy the conditions in Lemma \ref{lem-inhoestimate-Fos} and $q,\tilde{q} \geq 2$. If $l_1\in[r_\alpha,\infty],\tilde{l}_2\in[\tilde{r}_\beta,\infty]$,  we have $$\|\mathcal{T}F\|_{L_t^q(\R;\dot{B}_{l_1}^{\bar{s}})}\lesssim \|F\|_{L_t^{\tilde{q}'}(\R;\dot{B}_{\tilde{l}_2'}^s)}\ ,\  s=\bar{s}+\frac{1}{\tilde{q}'}+\frac{n}{\tilde{l}_2'}-\frac{1}{q}-\frac{n}{l_1}-2\ ,\  \mathcal{T}=\mathcal{G}\mbox{ or }\mathcal{H}\ .$$
\end{cor}

\begin{prf}
	By Bernstein inequality, we have $$\|P_j(\mathcal{T}F)\|_{L^{l_1}}\lesssim 2^{jn(\frac{1}{r_\alpha}-\frac{1}{l_1})}\|P_j(\mathcal{T}F)\|_{L^{r_\alpha}}$$and$$\|P_jF\|_{L^{\tilde{r}_\beta^{'}}} \lesssim 2^{jn(\frac{1}{\tilde{r}_\beta}-\frac{1}{\tilde{l}_2})}\|P_j(F)\|_{L^{\tilde{l}_2^{'}}} \ . $$ Also, by Lemma \ref{lem-inhoestimate-Fos} and scaling, $$\|P_j(\mathcal{G}F)\|_{L^qL^{r_\alpha}}\lesssim 2^{-j(\frac{1}{q}+\frac{n}{r_\alpha}+2-\frac{1}{\tilde{q}'}-\frac{n}{\tilde{r}_\beta^{'}})}\|P_j F\|_{L^{\tilde{q}'}L^{\tilde{r}_\beta^{'}}} \ .$$ Therefore, by Minkowski inequality, 
	\begin{eqnarray*}
		\|\mathcal{T}F\|_{L^q\dot{B}_{l_1}^{\bar{s}}}&\leq&\|2^{j\bar{s}}\|P_j(\mathcal{T}F)\|_{L^qL^{l_1}}\|_{l^2(\Z)}\\
		&\lesssim&\|2^{j(\bar{s}+\frac{n}{r_\alpha}-\frac{n}{l_1})}\|P_j(\mathcal{T}F)\|_{L^qL^{r_\alpha}}\|_{l^2(\Z)}\\
		&\lesssim&\|2^{j(\bar{s}-\frac{n}{l_1}-\frac{1}{q}-2+\frac{1}{\tilde{q}'}+\frac{n}{\tilde{r}_\beta^{'}})}\|P_jF\|_{L^{\tilde{q}'}L^{\tilde{r}_\beta^{'}}}\|_{l^2(\Z)}\\
		&\lesssim&\|2^{js}\|P_jF\|_{L^{\tilde{q}'}L^{\tilde{l}_2^{'}}}\|_{l^2(\Z)}\\
		&\leq&\|F\|_{L^{\tilde{q}'}\dot{B}_{\tilde{l}_2'}^s} \ .
	\end{eqnarray*}
\end{prf}

\subsection{Key linear estimate}
Equipped with the Strichartz estimates, we are able to present the following key estimate for linear wave equations.

To choose the appropriate acceptable pairs, we distinguish two cases for the power $p$:\\
{\rm (i)} $p_{H^1}< p<2$ and $n>6$: for fixed $\delta \in (0,\sigma)$, 
\begin{equation}\label{eq-nlowp}\tag{C1}
	\begin{gathered}
		\frac{1}{q}=\frac{1}{q_0}=\frac{n-3}{2(n-1-p)} \ , \ 
		\frac{1}{r_0}=\frac {n-4}{2n}+\frac{(n-3)(p-1)+n(2-p)\frac{n-1}{n-1-\delta}}{2n(n-1-p)}\ ,\\
		\frac{1}{r}=\frac{1}{n}(\frac{2}{p-1}-\frac{n-3}{2(n-1-p)}) \ , \ s_0=\frac{n}{2}-\frac{1}{q_0}-\frac{n}{r_0}\ ,
	\end{gathered}
	\nonumber
\end{equation}
{\rm (ii)} $p>p_{H^1}$ and $p\geq 2$:
\begin{equation}\label{eq-nhighp}\tag{C2}
	\frac{1}{q}=\frac{1}{q_0}=\frac{1}{p} \ , \ \frac{1}{r}=\frac{1}{n}(\frac{2}{p-1}-\frac{1}{p}) \ , \  \frac{1}{r_0}=\half-(1+\frac{1}{p})\frac{1}{n} \ , \ s_0=1\ .\nonumber
\end{equation}

\begin{lem}\label{lem-linearwave}
	Let $n\geq 4$, $p> p_{H^1}$ and $q, r, q_0, r_0, s_0$ satisfy \eqref{eq-nlowp} or \eqref{eq-nhighp}. Suppose $(q_1, r_1)$ are given by \beeq\label{eq-q1}\frac{1}{q_1}=\frac{p-1}{q}+\frac{1}{q_0} \ , \ \frac{1}{r_1}=\frac{p-1}{r}+\frac{1}{r_0}\ .\eneq
	For any $(f,g)\in \dot{\mathcal{H}}^{s_c}\times\dot{\mathcal{H}}^{s_c-1}$, $F\in L_t^{q_1}(\R;\dot{\mathcal{B}}_{r_1}^{s_c-s_0})$,
	there is a unique weak solution $u$, such that $$u\in C_t^1(\R;S')\ , \ (u,\pt u)\in L_t^{q}(\R;\dot{\mathcal{B}}_{r}^{0}\times\dot{\mathcal{B}}_{r}^{-1})\mbox{ and }u\in L_t^{q_0}(\R;\dot{\mathcal{B}}_{r_0}^{s_c-s_0})\ ,$$ satisfying the wave equation $\Box u=F$  with initial data $(f,g)$. Moreover,
	the solution satisfies the estimate
	\begin{equation}\label{eq-Ltxq}
		\|u\|_{L_t^q\dot{B}_{r}^0}+\|u\|_{L_t^{q_0}\dot{B}_{r_0}^{s_c-s_0}}+\|\pt u\|_{L_t^{q}\dot{B}_{r}^{-1}}\lesssim\|(f,g)\|_{s_c}+\|F\|_{L_t^{q_1}\dot{B}_{r_1}^{s_c-s_0}}\ ,
	\end{equation}
	where the implicit constant is independent of the choice of $f,g,F$.
	In addition, the selection of $q, r, q_0, r_0, s_0$ ensures $ r\geq r_0$.
\end{lem}
\begin{prf}		
	\textbf{Estimates}: At first, we observe that
	our choice of $(q,r)$ satisfies $s_c=n/2-1/q-n/r$,
	and
	the necessary scaling-invariance condition in 
	\eqref{eq-Fos} is ensured by \eqref{eq-q1} and our choices of $s_c$ and $s_0$. Meanwhile, it is sufficient to establish the estimate for $u$, because we can write $\pt u=S(t)(g, \Delta f)+\mathcal{H}[\sqrt{-\Delta}F](t)$, when $f, g\in C_c^\infty(\R^{n})$   and $F\in C_c^\infty(\R^{n+1})$. Noticing that $0<s_c-s_0<n/r_1$, the estimates for $\pt u$ is obtained by dense argument once all regularity indices in the estimate are lowered by one. 
	
	Let us begin with the case
	\eqref{eq-nlowp}, for which we divide the proof into two parts: the homogeneous part and the inhomogeneous part. 
	
	For the homogeneous part,
	by Lemma \ref{lem-hoestimate-Tao}, we only need to check that $(q,r),(q_0,r_0)$ are $\sigma$-admissible. 

	At first, as $p\in (p_{H^1},2)$, we notice that
	$q=q_0=2\frac{n-1-p}{n-3}>2$,
	$$\frac{1}{r_0}\ge
	\left.\frac{1}{r_0}\right|_{\de=0}
	=\frac{n-1}{2n}-\frac{n-3}{n(n-1-p)}\ ,$$
	$$			\frac{1}{r}\leq\frac{n-1}{2n}-\frac{n-3}{n(n-1-p)} \Leftrightarrow (n-1)(p-1)^2-(n^2-4n+9)(p-1)+4(n-2)\leq0\ .$$
	Denote the function $H(p):=(n-1)(p-1)^2-(n^2-4n+9)(p-1)+4(n-2)$. 
	As $n>6$, we observe that
	$$
	H(2)=-n^2+9n-18\leq 0,\ 
	H(p_{H^1})=\frac{4(6-n)}{(n-2)^2}\leq 0\ ,$$
	Thus we have $H(p)\le \max( H(p_{H^1}), H(2))\le 0$ for  any $p\in [p_{H^1},2)$. 
	In addition,
	as $\frac{1}{r}=\frac{1}{n}(\frac{2}{p-1}-\frac{n-3}{2(n-1-p)})$ is a decreasing function of $p>1$, it is clear that
	$$\frac{1}{r}\ge \left.\frac{1}{r}\right|_{p=2}=\frac{1}{n}(\frac{2}{2-1}-\frac{n-3}{2(n-1-2)})=\frac{3}{2n}>0\ .$$
	In particular,  $\frac{1}{r_0}\ge \frac{1}{r}>0$.
	
	Based on this fact,
	we need only to check that $1/r_0\le 1/2-2/((n-1)q_0)$.
	For this purpose, 
	we introduce the auxiliary function $r_\al$, which is defined by
	\beeq\frac{1}{r_\al}=\frac{n-5}{2(n-1)}+\frac{2(n-3)(p-1)+(n-1)(2-p)\frac{n-1}{n-1-\delta}}{2(n-1)(n-1-p)}\ ,\eneq
	and claim that \beeq\label{eqral}
	\frac{1}{r_0}\leq \frac{1}{r_\al}
	\leq \half-\frac{2}{(n-1)q_0}\ .\eneq
	The first inequality is equivalent to $p\geq 1+4(n-2)/(n-1)^2$, which is true for $p\geq p_{H^1}=1+4/(n-2)$. Noticing that $\delta<\sigma$ and $1/r_\al$ is increasing function of $\de$, we have
	\begin{equation}\label{eq-r_alpha}
		\begin{aligned}\frac{1}{r_\al}		<
			\left.\frac{1}{r_\al}\right|_{\de=\sigma}=
			&\frac{n-5}{2(n-1)}+\frac{2(n-3)(p-1)+2(n-1)(2-p)}{2(n-1)(n-1-p)}\\
			=&\half-\frac{n-3}{(n-1)(n-1-p)}=\half-\frac{2}{(n-1)q_0}\ .
		\end{aligned}
	\end{equation}

	Turning to the inhomogeneous part, for which we use Corollary \ref{cor-inhoesti-wave}. 
	Recalling that, as $p\in (p_{H^1},2), \delta \in (0,\sigma)$ and $n>6$, we have
	\begin{equation}
		\begin{gathered}
			\frac{1}{q_1}=\frac{p-1}{q}+\frac{1}{q_0}=\frac{(n-3)p}{2(n-1-p)}\in (\frac 12, 1)\ ,\\ 
			\frac{1}{r_1}=\frac{p-1}{r}+\frac{1}{r_0}= \half +\frac{2-p}{2(n-1-p)}\frac{n-1}{n-1-\delta} \in (\frac 12, 1)\ .
		\end{gathered}\nonumber
	\end{equation}
	In particular, together with \eqref{eq-r_alpha},  we see that
	$q_1, r_1\in (1,2)$ and $q, r_\alpha\in (2,\infty)$.
	
	To apply Corollary \ref{cor-inhoesti-wave}, we 
	introduce  the auxiliary pairs $(q, r_\al)$ and
	$(\tilde q,  \tilde r_\be) $\\$=(q_1', r_1')$. 
	Let $l_2=\tilde r_\be$ and $l_1=r, r_0$.
	Because the proof in the homogeneous part has shown $r\ge r_0\ge r_\al$, 
	to complete the proof for the inhomogeneous part, it remains to
	check that 
	the pairs $(q, r_\al)$ and $(\tilde q,  \tilde r_\be)$ satisfy the conditions of 
	Lemma \ref{lem-inhoestimate-Fos}. 
	
	Actually, by \eqref{eq-r_alpha},
	$(q, r_\al)$ is $\sigma$-admissible. For
	$(\tilde q,  \tilde r_\be)$,
	it is $\sigma$-acceptable:
	$$		\frac{1}{\tilde{q}}-2\sigma(\half-\frac{1}{\tilde{r}_\beta})
	=-\delta\frac{(n-1)(2-p)}{2(n-1-\delta)(n-1-p)}<0\ .
	$$
	It remains to check
	\eqref{eq-Ltxq}.
	At first, we have
	$$		\sigma(1-\frac{1}{r_\alpha}-\frac{1}{\tilde{r}_\beta})
	=\sigma(\frac{1}{r_1}-\frac{1}{r_\al})
	=1-\frac{(n-3)(p-1)}{2(n-1-p)}=1-\frac{p-1}q=\frac{1}{q}+\frac{1}{\tilde{q}}<1\ . 
	$$
			Similar to $1/r_\al$,
			$1/r_1$ is also increasing with respect to $\de$, then we have
			$$\frac{1}{\tilde{r}_\be}=1-\frac{1}{r_1}\le
			1-\left.\frac{1}{r_1}\right|_{\de=0}
			=\half -\frac{2-p}{2(n-1-p)}=
			\frac{n-3}{2(n-1-p)}\ ,
			$$
			$$
			\frac{1}{r_\al}		\ge
			\left.\frac{1}{r_\al}\right|_{\de=0}
			=\frac{n-5}{2(n-1)}+\frac{(n-5)(p-1)+n-1}{2(n-1)(n-1-p)}
			=\frac{(n-3)^2}{2(n-1)(n-1-p)}
			\ .
			$$
			Thus $$\frac{r_\alpha}{\tilde{r}_\be}\leq \left.\frac{r_\alpha}{\tilde{r}_\be}\right|_{\delta=0}
			=\frac{n-1}{n-3}
			=\frac{\sigma}{\sigma-1}\ .$$ 
			Since $$\frac{1}{\tilde{q}}=1-\frac{p}q \leq\frac{1}{q}\Leftrightarrow p\geq 1+\frac{2}{n-1}\ , $$ 
			$$\frac{1}{\tilde{q}}=1-\frac{p}q
			=\frac{n-1}{2}\frac{2 - p}{n-1-p}=
			\sigma
			(\left.\frac{1}{r_1}\right|_{\de=\sigma}-\frac 12)
			\ge \sigma
			(\frac{1}{r_1}-\frac 12)
			=\sigma(\frac 12- \frac{1}{\tilde r_\be} )\ ,
			$$
			we obtain that 
			$$\frac{1}{r_\alpha}\leq \half-\frac{1}{\sigma q}\leq\half-\frac{1}{\sigma \tilde{q}}
			\leq \frac{1}{\tilde{r}_\be} \Rightarrow \frac{r_\alpha}{\tilde{r}_\be}\ge 1\ge \frac{\sigma-1}{\sigma}\ , $$ which verifies \eqref{eq-Fos}.
			This completes the proof of  \eqref{eq-Ltxq}   in the case of  \eqref{eq-nlowp}.

			Next, we consider the Case \eqref{eq-nhighp}. 
			In this case, we will apply  Lemma \ref{lem-hoestimate-Tao}, for which it suffices to   show that $(q,r)$, $(q_0,r_0)$ and $(q_1',r_1')$ are $\sigma$-admissible.
			
			As $q=q_0=p\ge 2$, we claim $$\frac{1}{p}\leq\sigma(\half-\frac{1}{r_0})\leq\sigma(\half-\frac{1}{r})\ .$$
			For the first inequality, $$\frac{1}{p}\leq\sigma(1+\frac{1}{p})\frac{1}{n} \Leftrightarrow p\geq \frac{n+1}{n-1}\ ,$$ which is true for $p\geq \frac{n+2}{n-2}=p_{H^1}$.
			The second inequality is equivalent to
			$$
			\frac{1}{r_0}\geq\frac{1}{r}
			\Leftrightarrow
			p\geq \frac{n+2}{n-2}\ .
			$$

			At last, in view of \eqref{eq-q1} and \eqref{eq-nhighp}, we have
			$$
			\frac{1}{q_1'}=1-\frac 1{q_1}=1-\frac{p-1}q-\frac 1{q_0}=0 \ ,  \ \frac{1}{r_1'}=1-\frac{p-1}r-\frac 1{r_0}=\half\ , 
			$$ which shows that $(q_1',r_1')=(\infty, 2)$ is $\sigma$-admissible.
			This completes the proof.
			
			\textbf{$C^1$ regularity}: In both two cases, we have $-n/r'<-2<0<n/r$, $-n/r_1'<s_c-s_0<n/r_1$ and $-n/2<s_c-1<s_c<n/2$. By Remark \ref{rem-homoBreal} and Proposition \ref{prop-density}, $C_c^\infty(\R^n)$ is dense in $ \dot{\mathcal{H}}^{s_c}$, $\dot{\mathcal{H}}^{s_c-1}$, $\dot{\mathcal{B}}_{r'}^{1}$, $\dot{\mathcal{B}}_{r'}^{2}$, $\dot{\mathcal{B}}_{r_1'}^{s_0-s_c}$, and $C_c^\infty(\R^{n+1})$ dense in $L_t^{q_1}(\dot{\mathcal{B}}_{r_1}^{s_c-s_0})$.
			We choose sequences of functions, $f_n, g_n \in C_c^\infty(\R^n)$  and $F_n \in C_c^\infty(\R^{n+1})$ such that $\lim f_n=f$, $\lim g_n=g$ and $\lim F_n=F$. Denote $u_n$ the solution to the wave equation $\Box u_n=F_n$  with initial data $(f_n,g_n)$.
			Then, we have, for all $\phi \in S(\R^n)$,
			
			$$\<u_n(t_2),\phi\>-\<u_n(t_1),\phi\>=\int_{t_1}^{t_2}\<\pt u_n(\tau),\phi\>d\tau\ ,$$
			\begin{equation}\label{eq-u-eqcontinuity}
				|\<u_n(t_2),\phi\>-\<u_n(t_1),\phi\>|\leq \|\pt u_n\|_{L_t^{q}\dot{B}_{r}^{-1}} \|\phi\|_{\dot{B}_{r'}^{1}}|t_2-t_1|^{1\over q'}\ ,
			\end{equation}
			$$|\<u_n(t),\phi\>|\leq \|\pt u_n\|_{L_t^{q}\dot{B}_{r}^{-1}} \|\phi\|_{\dot{B}_{r'}^{1}}|t|^{1\over q'}+\|f_n\|_{\dot{H}^{s_c}}\|\phi\|_{\dot{H}^{-s_c}}\ . $$
			Hence for all time $t$, $u_n(t)$ has a limit $\bar{u}(t)$ in $S'$ and $\bar{u}\in C(S')$. Also, $\bar{u}$ can be identified as the distribution, $$\<\bar{u},\psi\>=\int_{\R}\<\bar{u}(t),\psi(t)\>_xdt\ , \forall \psi \in C_c^\infty(\R^{n+1})\ .$$ So, on the one hand, for all $\psi \in C_c^\infty(R^{n+1})$ we have $$|\<\bar{u},\psi\>-\<u_n,\psi\>|\leq C_\psi \|\pt u_n-\pt u\|_{L_t^{q}\dot{B}_{r}^{-1}}\rightarrow 0\ .$$ On the other hand, Strichartz estimates tell us that $u_n\rightarrow u$ in $\mathcal{D'}(\R^{n+1})$, the dual space of $C_c^\infty(\R^{n+1})$. Thus $u=\bar{u}\in C(S')$. 
			
			The argument for the continuity of $\pt u$ is almost the same. Noticing that $\partial_{tt}^2u_n=\Delta u_n+ F_n$ and $u_n$ is a Cauchy sequence in $L_t^{q}\dot{\mathcal{B}}_{r}^{0}$, we have some similar inequalities, 
			\begin{align}
				|\<\pt u_n(t_2)-\pt u_n(t_1),\phi\>|  \leq &\|u_n\|_{L_t^{q}\dot{B}_{r}^{0}}
				\|\phi\|_{\dot{B}_{r'}^{2}}|t_2-t_1|^{1\over q'}		 \nonumber
				\\
				+\| F_n\|_{L_t^{q_1}([t_2,t_1])\dot{B}_{r_1}^{s_c-s_0}} &\|\phi\|_{\dot{B}_{r_1'}^{s_0-s_c}}|t_2-t_1|^{1\over q_1'},\ 					\forall \phi \in S(\R^n)\ .
				\label{eq-ptu-eqcontinuity}			
			\end{align}
			Therefore, by the same reasoning, we obtain $\pt u \in C(S')$.
			
			\textbf{Uniqueness}: Assume there is a solution $u\in C^1(S')$ to the wave equation such that $f,g=0$ and $F=0$. Then, because $\p_t^{k+2} u= \Delta \p_t^ku$, $u$ belongs to $C^\infty(S')$, hence the uniqueness.
		\end{prf}
		\begin{rem}\label{rem-s0}
			When $p<2$,
			the index $s_0$ in \eqref{eq-nlowp} is $$s_0=2-\frac{(n-3)p+n(2-p)\frac{n-1}{n-1-\delta}}{2(n-1-p)}\ .$$ As we could set $\delta\in(0,\sigma)$ arbitrarily small, the limit of our $s_0$ is
			$$\bar{s}_0=
			2-\frac{(n-3)p+n(2-p)}{2(n-1-p)}
			=	2-\frac{2n-3p}{2(n-1-p)}
			=\frac 12+\frac{n-3}{2(n-1-p)}	
			\ ,$$
			which is the regularity index \eqref{eq-s0}, appeared in Theorem \ref{thm-main}.
		\end{rem}
		\begin{rem}
			The $C^1$ regularity ensures that for any solution $u$ in Lemma \ref{lem-linearwave}, $$\<u \mathbb{E}_+, \Box\psi\>=-\<f,\pt\psi(0)\>_x+\<g,\psi(0)\>_x+\<F \mathbb{E}_+,\psi\>\ , \forall \psi \in C_c^\infty(\R^{n+1})\ , $$
			where $\mathbb{E}_+$ is the characteristic function of the time interval $(0,\infty)$.
		\end{rem}
		
		\subsection{Nonlinear estimates}
		\begin{lem}[\cite{MR3468325}, Proposition 2.1.2]\label{lem-fraLei}
			Let $p\in (1,\infty)$ and $0< s< n/p$. We have the following fractional Leibniz rule:
			$$\|uv\|_{\dot{B}_p^s}\lesssim(\|v\|_{L^\infty}+\|v\|_{\dot{B}_{p,\infty}^{n/p}})\|u\|_{\dot{B}_p^s}\ .$$ 
		\end{lem}
		
		Now we infer that solutions in Lemma \ref{lem-linearwave} can be localized. This trick is also used in \cite{LdSo95}. Denote the set $\Lambda_{R,0}:=\{(t,x)\in\R^{1+n}| \ |x|<R-t, t>0\}$, $0<R< \infty$.
		
		\begin{lem}\label{lem-localised}
			Let $\chi$ be a cut-off function in $\R^{n+1}$ and equal to $1$ in a neighborhood of $\Lambda_{R,0}$. If $w$ and $w_\chi$ are the solutions in Lemma \ref{lem-linearwave} satisfying $\Box w= F$, $\Box w_\chi= \chi F$ with zero initial data, then $w=w_\chi$  in $\Lambda_{R,0}$.
		\end{lem}
		
		\begin{prf} By Lemma \ref{lem-fraLei}, we have $\chi F\in L_t^{q_1}\dot{\mathcal{B}}_{r_1}^{s_c-s_0}$. So for all $\psi \in 
			C_c^\infty(\R^{n+1})$, the following equality is valid $$\int_{\R\times \R^n} (w-w_\chi)\mathbb{E}_+ \Box\psi \ dtdx=\int_{\R\times \R^n} (1-\chi)F \mathbb{E}_+\psi \ dtdx\ .$$
			Because for all $\phi\in C_c^\infty(\Lambda_{R,0})$, one can find some $\psi\in C_c^\infty{(\R^{n+1})}$ such that $\Box \psi=\phi$ in $\{(t,x)|t>0\}$ and $(1-\chi) \psi=0$. Then $$\int_{\R\times \R^n} (w-w_\chi) \mathbb{E}_+ \Box\psi \ dtdx=\int_{\R\times \R^n} (w-w_\chi) \mathbb{E}_+\phi \ dtdx=0 \ ,\forall \phi \in C_c^\infty(\Lambda_{R,0})\ ,$$ which shows that
			$w=w_\chi$  in $\Lambda_{R,0}$.
		\end{prf}

		\begin{lem}[\cite{NO97}, Lemma 2.2]\label{lem-fracdir-Besov} Assume $1< p < \infty$ and $0<s<p$. Let $r,m,l$ satisfy $1<r\leq 2 \leq l$, $m<\infty$, $1/r=(p-1)/m+1/l$. 
			We have the following fractional chain rule:
			\begin{equation}\label{eq-fracdir-Besov}
				\||u|^p\|_{\dot{B}_{r}^{s}}\lesssim \|u\|_{\dot{B}_{m}^{0}}^{p-1}\|u\|_{\dot{B}_{l}^{s}}\ .
			\end{equation}
		\end{lem}
		
		\begin{cor}\label{cor-F(u)}
			Let $p>p_{H^1}$ and $r, r_0, r_1, s_0$ be given in \eqref{eq-nlowp}-\eqref{eq-nhighp} and \eqref{eq-q1}. Assume that $s_c-s_0\in (0,p)$ is fulfilled.
			Then,
			for any $u\in \dot{\mathcal{B}}_r^0\cap\dot{\mathcal{B}}_{r_0}^{s_c-s_0}$, we have
			$|u|^p\in \dot{\mathcal{B}}_{r_1}^{s_c-s_0}$ which enjoys the
			estimate \eqref{eq-fracdir-Besov}.
			
		\end{cor}
		\begin{prf} 
			It is obvious that $1<r_1\leq 2 \leq r,r_0<\infty$. So, by Lemma \ref{lem-fracdir-Besov}, we only need to verify that $|u|^p\in S_h'$. Noticing that $\dot{\mathcal{B}}_r^0\hookrightarrow L_r$, we have $|u|^p\in L_{r/p}$, 
			$r/p\geq\frac{(p-1)(n+1)}{2p}>1$, which means $|u|^p \in S_h'$.\end{prf}

		\section{Proof of Theorem \ref{thm-main}}
		
		Actually, we can prove the following strengthened version of Theorem \ref{thm-main}.
		\begin{thm}\label{thm-mainv2}
			Let $n\ge 4$, $p>p_{H^1}$,
			$s_c= n/ 2- {2}/(p-1)$ and $\bar s_0$ be given in
			\eqref{eq-s0}.
			Suppose we have $p>s_c-\bar s_0$,
			for any given $(f, g)\in \dot{\mathcal{H}}^{s_c}\times \dot{\mathcal{H}}^{s_c-1}$, there exists $T>0$ such that there exists a unique weak solution $u\in L^q_tL^r_x$ to
			$\Box u=|u|^p$ in $(0,T)\times\R^n$ with data $u(0)=f$, $u_t(0)=g$,
			where $(q,r)$ is given in
			\eqref{eq-nlowp}-\eqref{eq-nhighp} (with $\delta\ll 1$ such that $p>s_c-s_0$).
			Moreover, the solution satisfies
			$u\in
			L_t^q([0,T);\dot{\mathcal{B}}_r^0)\cap L_t^{q_0}([0,T);\dot{\mathcal{B}}_{r_0}^{s_c-s_0})$ and $u\in C_t^1([0,T);S')$.
			In addition, the solution is global if the data is small in $\dot{\mathcal{H}}^{s_c}\times \dot{\mathcal{H}}^{s_c-1}$.
		\end{thm}

		\subsection{Existence} 
		As we have noticed in
		Remark \ref{rem-s0}, when $p<2$, $\lim_{\de\to 0}s_0=\bar s_0$ and so there exists $\de\in (0,\sigma)$ such that
		$p>s_c-s_0$. In the following, we fix the choice of $\de$ for the case $p<2$.
		
		Let $(q, r)$, $(q_0, r_0)$ and $(q_1, r_1) $ be given in \eqref{eq-nlowp}-\eqref{eq-nhighp} and \eqref{eq-q1}. 
		Based on Lemma \ref{lem-linearwave} and Corollary \ref{cor-F(u)},
		we will use a standard iteration argument to give the proof.

		Given $(u(0) , \partial_t u(0))=(f, g)\in \dot{\mathcal{H}}^{s_c}\times \dot{\mathcal{H}}^{s_c-1}$, with $T>0$ to be determined, we introduce a Banach space
		$$X_T:=L_t^q\dot{\mathcal{B}}_r^0\cap L_t^{q_0}\dot{\mathcal{B}}_{r_0}^{s_c-s_0}([0,T)\times\R^n)
		\ .$$
		Let $$\Phi(u)(t,x)=\mathbb{E}_{[0,T)}(t)\times(
		S(t)(f, g)+\mathcal{G}(\mathbb{E}_{[0,T)}|u|^p)(t))
		\ ,$$where $\mathbb{E}_{[0,T)}$ is the characteristic function of the time interval $[0,T)$.
		By Lemma \ref{lem-linearwave}, Lemma \ref{lem-fracdir-Besov} and H\"older's inequality, there exists $C>0$, such that for any $u\in X_T$, we have
		\begin{eqnarray}
			\|\Phi(u)(t)\|_{X_T}&\le&
			\|S(t)(f, g)\|_{X_T}+C\|u\|_{L_t^q\dot{B}_{r}^{0}}^{p-1}\|u\|_{L_t^{q_0}\dot{B}_{r_0}^{s_c-s_0}} \nonumber\\
			&	\le &
			\|S(t)(f, g)\|_{X_T}+C\|u\|_{X_T}^p\ .\label{eq-main-iter}
		\end{eqnarray}
		
		Let $\ep\le \ep_0=(2^p C)^{-1/(p-1)}$, be such that $C(2\ep )^p\le \ep$.
		By Lemma \ref{lem-linearwave}, 
		$$\|S(t)(f, g)\|_{X_\infty}\les \|(f,g)\|_{s_c}\ .$$
		Due to the fact that $q,q_0<\infty$, 
		we see that there exists $T(\ep)>0$, such that 
		\beeq\label{eq-initial} \|S(t)(f, g)\|_{X_T}\le \ep\ ,\eneq
		for any $T\in (0, T(\ep)]$.
		In particular, we could choose $T(\ep)=\infty$, if $\|(f,g)\|_{s_c}$ is sufficiently small.
		
		We take $u_{-1}\equiv 0$ and define iteratively $u_m(t)=\Phi(u_{m-1})(t)$, $m\in \N$. By \eqref{eq-main-iter} and
		\eqref{eq-initial}, for any $m\in\N$, we have $$u_m\in X_T, \ \|u_m\|_{X_T}\leq 2\ep   \ .$$
		Then
		$|u_m|^p\in L_t^{q_1}\dot{\mathcal{B}}_{r_1}^{s_c-s_0}$ by Corollary \ref{cor-F(u)}.
		
		Noticing that $||z_1|^p-|z_2|^p|\lesssim (|z_1|^{p-1} + |z_2|^{p-1})|z_1-z_2|$ and all $u_m$ can be localised, by Lemma \ref{lem-localised}, then for any fixed $R$ and any cut-off function $\chi$ with $\chi$ equalling to $1$ in a neighborhood of $\Lambda^T_{R,0}:=\Lambda_{R,0}\cap ([0,T)\times\R^n)$, we have, for some $C_1>1$
		\begin{eqnarray}
			&&\|u_m-u_k\|_{L_t^{q_0}L_x^{r_0}(\Lambda^T_{R,0})}\nonumber\\
			&=& \|(u_m-u_k)_\chi\|_{L_t^{q_0}L_x^{r_0}(\Lambda^T_{R,0})}\nonumber\\
			&\leq& \|(u_m-u_k)_\chi\|_{L_t^{q_0}L_x^{r_0}([0,T)\times\R^n)}\nonumber\\
			&\le &C_1\|\chi (F(u_{m-1})-F(u_{k-1}))\|_{L_t^{q_1}L_x^{r_1}([0,T)\times\R^n)}\nonumber\\
			&\le &C_1^2 (\|u_{m-1}\|_{L_t^{q}L_x^{r}}^{p-1}+\|u_{k-1}\|_{L_t^{q}L_x^{r}}^{p-1})\|\chi (u_{m-1}-u_{k-1})\|_{L_t^{q_0}L_x^{r_0}}\nonumber\\
			&\le &C_1^3 (\|u_{m-1}\|_{L_t^q\dot{B}_{r}^0}^{p-1}+\|u_{k-1}\|_{L_t^q\dot{B}_{r}^0}^{p-1})\|\chi(u_{m-1}-u_{k-1})\|_{L_t^{q_0}L_x^{r_0}}\nonumber\\
			&\leq& \half\|\chi(u_{m-1}-u_{k-1})\|_{L_t^{q_0}L_x^{r_0}}\ ,\label{eq-iteration}
		\end{eqnarray}
		provided that $2C_1^3 (2\ep)^{p-1}\le 1/2$, i.e. $\ep\le \ep_1>0$ with
		$2C_1^3 (2\ep_1)^{p-1}= 1/2$.
		As $C_1$ is independent of $\chi$, by shrinking $\supp \chi$ to $\Lambda_{R,0}$,  we conclude that
		$$
		\|u_m-u_k\|_{L_t^{q_0}L_x^{r_0}(\Lambda^T_{R,0})}\leq \frac{1}{2}\|u_{m-1}-u_{k-1}\|_{L_t^{q_0}L_x^{r_0}(\Lambda^T_{R,0})}\ .$$
		Also, since $r_0\leq r$ and $q=q_0$, we have
		\begin{equation}
			\|u_0\|_{L_t^{q_0}L_x^{r_0}(\Lambda^T_{R,0})}\lesssim_R \|u_0\|_{L_t^qL_x^r(\Lambda^T_{R,0})}<\infty\ .\nonumber
		\end{equation}
		Therefore, there exists some $u\in  L_t^{q_0}L_{x, loc}^{r_0}$ such that $u_m\rightarrow u$ in   $L_t^{q_0}L_x^{r_0}$ locally, thus in $\mathcal{D'}$. As mentioned in Remark \ref{rem-homoBreal}
		, $\dot{\mathcal{B}}_r^0$ and $\dot{ \mathcal{B} }_{r_0}^{s_c-s_0}$ are reflexive Banach spaces, due to the fact that $-n/r_0'<(s_c-s_0)<n/r_0$. Since $\{u_m\}$ is bounded in both $L_t^q\dot{\mathcal{B}}_r^0$ and $L_t^{q_0}\dot{\mathcal{B} }_{r_0}^{s_c-s_0}$, there exists a subsequence $\{u_{m_k}\}$ such that,
		$$u_{m_k} \rightharpoonup  \tilde{u}\  \mbox{ in } L_t^q\dot{\mathcal{B}}_r^0
		\cap L_t^{q_0}\dot{\mathcal{B} }_{r_0}^{s_c-s_0}
		\Longrightarrow u_{m_k}\rightarrow \tilde{u} \mbox{ in } \mathcal{D'}\ ,$$
		This means $u=\tilde{u} \in X_T$. Thus $F(u)\in L_t^{q_1}\dot{\mathcal{B}}_{r_1}^{s_c-s_0}$ and $u\in L_t^qL_x^r$, by Corollary \ref{cor-F(u)} and embedding.
		
		Similarly, by repeating the argument from the third line in \eqref{eq-iteration}, we have $F(u_m)\rightarrow F(u)$ in $L_t^{q_1}L_x^{r_1}$  locally. Therefore $u$ is the desired weak solution in $[0,T)\times \R^n$, provided $\ep\le \min(\ep_0,\ep_1)$.
		\subsection{$C^1$ regularity} In the case \eqref{eq-nlowp}, the inequalities \eqref{eq-u-eqcontinuity} and \eqref{eq-ptu-eqcontinuity} show that for all $\phi \in S$ the function sequences $\{\<u_n(t),\phi\>_x\}$  and $\{\<\pt u_n(t),\phi\>_x\}$ are uniform bounded and equicontinuous in any finite time interval. Thus, by Arzel\`a–Ascoli theorem , there is a subsequence convergent to distributions $v, \pt v\in C(S')$. Because $u_n \rightarrow u$  in $\mathcal{D'}$, we have $u=v \in C_t^1(S')$. 
		
		In the case \eqref{eq-nhighp}, 
		by the inequality \eqref{eq-u-eqcontinuity}, we have $u\in C(S')$. But this time $1/q_1'$ is zero in \eqref{eq-ptu-eqcontinuity}. Notice that, for all $\phi\in C_c^\infty(\R^n)$ and $0\leq t_1<t_2<T$, 
		\begin{eqnarray*}
			\<\pt u_n(t_2)-\pt u_n(t_1),\phi\>_x=\<\Delta u_n,\mathbb{E}_{[t_1,t_2]}\phi\>_{t,x}+\< F(u_{n-1}),\mathbb{E}_{[t_1,t_2]}\phi\>_{t,x}\ ,\\
			\<\pt u_n(t_2),\phi\>_x=\<\Delta u_n,\mathbb{E}_{[0,t_2]}\phi\>_{t,x}+\< F(u_{n-1}),\mathbb{E}_{[0,t_2]}\phi\>_{t,x}+\<g,\phi\>_x\ . 
		\end{eqnarray*} The proof in the existence part shows that $u_m\rightarrow u$ in $L_t^{q_0}L_x^{r_0}$ locally and $F(u_m)\rightarrow F(u)$ in $L_t^{q_1}L_x^{r_1}$  locally. Thus, there exists a distribution $w\in C([0,T);\mathcal{D'})$ such that for all fixed $t\in [0,T), \pt u_n(t)\rightarrow w(t)$ in $\mathcal{D'}$. Moreover, $\pt u_n\rightarrow w$ in $\mathcal{D'}((0,T)\times \R^n)$, hence $\pt u= w$ in $(0,T)\times \R^n$. Finally, according to the energy estimate, $\pt u$ belongs to $L^\infty(\dot{H}^{s_c-1})$, which infers that $\pt u \in C(S')$.
		
		\subsection{Uniqueness} Let $S_T=[0,T)\times\R^n$ and $u_1, u_2\in L_t^qL_x^r$ be two solutions to the equation \eqref{eq-Strauss}.
		Then $u\in  L_t^{q_0}L_{x, loc}^{r_0}$ and
		we have
		$$
		\|u_1-u_2\|_{L_t^{q_0}L_x^{r_0}(\Lambda^T_{R,0})}\leq C(\|u_1\|_{L_t^qL_x^r(S_T)}^{p-1}+\|u_2\|_{L_t^qL_x^r(S_T)}^{p-1}) \|u_1-u_2\|_{L_t^{q_0}L_x^{r_0}(\Lambda^T_{R,0})}$$
		for all $R\in(0,\infty)$.
		By choosing $T\ll 1$ such that $$C(\|u_1\|_{L_t^qL_x^r(S_T)}^{p-1}+\|u_2\|_{L_t^qL_x^r(S_T)}^{p-1})\leq 1/2 \ ,$$ we deduce that $\|u_1-u_2\|_{L_t^{q_0}L_x^{r_0}(\Lambda^T_{R,0})}=0$. By spatial translation invariance,  we obtain $\|u_1-u_2\|_{L_t^{q_0}L_x^{r_0}(S_T)}=0$.
		This proves the uniqueness for general $T>0$, by
		iterating the above process.
		
		\begin{rem}[Sharpness of the $\bar{s}_0$]
			In \eqref{eq-nhighp}, 
			$\bar{s}_0=1$ and it is sharp,
			in view of the wellposedness in $H^{s_c}$. Here, in the case of \eqref{eq-nlowp}, we try to explain the sharpness of our choice
			of $\bar{s}_0$,
			in the framework of 
			\eqref{eq-Ltxq},
			\eqref{eq-q1},
			and 
			Foschi's inhomogeneous Strichartz estimates  \cite{MR2134950}.
			Actually, the sharpness of the $\bar{s}_0$ is equivalent to the following linear programming problem:$$\sup \frac{1}{q_1'}+\frac{n}{r_1'}$$
			\begin{equation}\label{eq-LNfull}
				\left\{
				\begin{gathered}
					\frac{1}{r_\alpha}\geq \frac{1}{r}\  , \ \frac{1}{r_\beta}\geq \frac{1}{r_0}\ ,\\ \frac{1}{r_\gamma'},\frac{1}{r_\zeta'}\geq\frac{1}{r_1'}\ ,\\
					\frac{1}{q}+\frac{1}{q_1'}=\frac{n-1}{2}(1-\frac{1}{r_\alpha}-\frac{1}{r_\gamma'})<1\ ,\  
					\frac{1}{q_0}+\frac{1}{q_1'}=\frac{n-1}{2}(1-\frac{1}{r_\beta}-\frac{1}{r_\zeta'})<1\ ,\\
					0\leq\frac{1}{q}\leq\min(\frac{n-1}{2}(\half-\frac{1}{r_\alpha}),\half)\ , \ 0\leq\frac{1}{q_0}\leq\min(\frac{n-1}{2}(\half-\frac{1}{r_\beta}),\half)\ ,\\ 0<\frac{1}{q_1'}<\min((n-1)(\half-\frac{1}{r_\gamma'}),\half)\ , \ 0<\frac{1}{q_1'}<\min((n-1)(\half-\frac{1}{r_\zeta'}),\half)\ ,\\
					\frac{n-3}{n-1}\leq\frac{r_\alpha}{r_\gamma'}\leq\frac{n-1}{n-3}\ ,\ \frac{n-3}{n-1}\leq\frac{r_\beta}{r_\zeta'}\leq\frac{n-1}{n-3}\ ,\\
					1-\frac{p-1}{q}=\frac{1}{q_0}+\frac{1}{q_1'}\ , \ 1-\frac{p-1}{r}=\frac{1}{r_0}+\frac{1}{r_1'}\ ,\ \frac{2}{p-1}=\frac{1}{q}+\frac{n}{r}\ .
				\end{gathered}
				\right.
			\end{equation}
			This is a little bit complicated. Fortunately, if we delete some subjections and consider a simplified linear programming problem, $$\sup \frac{1}{q_1'}+\frac{n}{r_1'}$$
			\begin{equation}\label{eq-LNsimple}
				\left\{
				\begin{gathered}
					\frac{1}{r_\gamma'},\frac{1}{r_\zeta'}\geq\frac{1}{r_1'}>0 \ ,\ \frac{1}{q}\geq0\ ,\\
					\frac{1}{q}+\frac{1}{q_1'}=\frac{n-1}{2}(1-\frac{1}{r_\alpha}-\frac{1}{r_\gamma'})\ ,\ 
					\frac{1}{q_0}+\frac{1}{q_1'}=\frac{n-1}{2}(1-\frac{1}{r_\beta}-\frac{1}{r_\zeta'})\ ,\\ 0<\frac{1}{q_1'}<(n-1)(\half-\frac{1}{r_\zeta'})\ ,\\
					\frac{n-3}{n-1}\leq\frac{r_\alpha}{r_\gamma'}\leq\frac{n-1}{n-3}\ ,\ \frac{n-3}{n-1}\leq\frac{r_\beta}{r_\zeta'}\leq\frac{n-1}{n-3}\ ,\\
					1-\frac{p-1}{q}=\frac{1}{q_0}+\frac{1}{q_1'}\ , \ 1-\frac{p-1}{r}=\frac{1}{r_0}+\frac{1}{r_1'}\ ,\ \frac{2}{p-1}=\frac{1}{q}+\frac{n}{r}\ ,
				\end{gathered}
				\right.
			\end{equation}
			we will find the answer is also $\bar{s}_0$. Our proof in Lemma \ref{lem-linearwave} shows that $\bar{s}_0$ is the limit of the $s_0$. Noticing the supremum in the linear programming \eqref{eq-LNfull} is no more than in \eqref{eq-LNsimple}, we deduce that the supremum in \eqref{eq-LNfull} is $\bar{s}_0$.
		\end{rem}
		
		\section{Appendix: Function Spaces}
		\subsection{Homogeneous Besov spaces}
		We introduce the realization of the homogeneous Besov space  $\dot{\mathcal{B}}_{p,q}^s$, a subspace of the tempered distribution $S'(\R^n)$, see, e.g.,
		\cite{BCD} and \cite{Tri}. Though there are many realizations of homogeneous Besov spaces. When $s<n/p$,  $\dot{B}_{p,q}^s$ admits a unique dilation commuting realization in $S'$ \cite[Theorem 4.1]{Bou}.
		
		\begin{defn}[\cite{BCD}]\label{defn-realispace}
			We denote by $S_h^{'}(\R^n)$ the space of all $u\in S'(\R^n)$ such that $$\lim_{\lambda \to \infty} \|\theta(\lambda \sqrt{-\Delta})u\|_{L^\infty}=0 \mbox{ for some } \theta \in C_c^{\infty}(R^n) \mbox{ with }\theta(0)\neq 0\ .$$
		\end{defn}

		\begin{rem}\label{rem-realispace}
			If for some $\theta\in C_c^\infty$, $\theta(0)\neq 0$, we have $\theta(\sqrt{-\Delta}))u \in L^p$, $p\in[1,\infty)$, then $u$ belongs to $ S_h'$. Also, any nonzero polynomial does not belong to $S_h^{'}$. See \cite[Definition 1.26, p.22]{BCD}.\end{rem}
		
		\begin{defn}\label{defn-homoBreal}
			Let $s\in \R$, $p\in[1,\infty]$, the realization of homogeneous Besov space $\dot{\mathcal{B}}_{p}^s$ consists of all the distributions $u \in S_h'(\R^n)$ such that
			$$\|u\|_{\dot{B}_{p}^s}=\|2^{js}\|P_ju\|_{L^p}\|_{l^2(\Z)}<\infty\ ,$$
			where $\{P_j\}_{j\in \Z}$ is the homogeneous Littlewood-Paley decomposition.
		\end{defn}
		\begin{rem}\label{rem-homoBreal}
			\begin{enumerate}
				\item When $p=2$, $\dot{\mathcal{B}}_2^s$ is the realization of homogeneous Sobolev space, denoted by $\dot{\mathcal{H}}^s$.
				\item The space $S_0$ of functions in $S$ whose Fourier transforms are supported away from $0$ is dense in $\dot{\mathcal{B}}_{p}^s(1\leq p<\infty)$. See \cite[Theorem 2.27, p.69]{BCD}.
				\item For $u\in S_h'$, $\phi \in S$ or $u\in S'$, $\phi \in S_0$, observe that $$\<u,\phi\>=\sum_{|j-j'|\leq1}\<P_ju,P_{j'}\phi\> \ .$$ Thus,  $|\<u,\phi\>|\lesssim\|u\|_{\dot{B}_p^s}\|\phi\|_{\dot{B}_{p'}^{-s}}$, for all $1\leq p\leq \infty$ and $s\in \R$. However, Schwartz function space is not always continuously embedded in $\dot{\mathcal{B}}_{p'}^{-s}$.
				\item  The realization of homogeneous Besov space $\dot{\mathcal{B}}_{p}^s$ is a Banach space if $(s,p)$ satisfies$$s<\frac{n}{p}\ ,\ p>1\ .$$ See \cite[Theorem 2.25, p.67]{BCD}.
				\item  If $1<p<\infty$, $-n/p'<s<n/p$, then $S\hookrightarrow\dot{\mathcal{B}}_p^s \hookrightarrow S'$ and, moreover, $S$, $C_c^\infty$ are dense in $\dot{\mathcal{B}}_p^s$. See \cite[Theorem 3.24]{Tri}.
				\item  If  $1<p<\infty$, $-n/p'<s<n/p$, then $\dot{\mathcal{B}}_{p}^s$ is reflexive.
				\item  If $2\leq p<\infty$, then $\dot{\mathcal{B}}_p^0\hookrightarrow L_p$.
			\end{enumerate}
		\end{rem}
		\begin{prf}
			We just show a sketch of proof for (6) and (7) in Remark \ref{rem-homoBreal}.
			
			(6): Let $T\in(\dot{\mathcal{B}}_p^s)'$. We have $|T\phi|\leq \|T\|\|\phi\|_{\dot{\mathcal{B}}_p^s}$, $\forall\phi \in \dot{\mathcal{B}}_p^s$.
			Meanwhile, as $S$ is continuously embedded in $\dot{\mathcal{B}}_p^s$, there exists a tempered distribution $u_T=l(T)$ such that $\<u_T,\phi\>=T\phi$, $\forall \phi \in S$. Let $u_N=\sum_{j\ge -N} P_j u_T\in S_h'$. Then, by Proposition 2.29 in \cite{BCD}, we have
			\begin{equation}\label{eq-eqnorm}
				\|(2^{-js}\|P_ju_N\|_{L^{p'}})\|_{l^2}\leq C_1 \sup_{\phi \in S, \|\phi\|_{\dot{\mathcal{B}}_p^s}\leq 1} |\<u_N,\phi\>|\le C_2 \|T\| \ ,
			\end{equation}
			for some uniform constants $C_1$, $C_2$. As $-s<n/p'$, \cite[Remark 2.24, p.66]{BCD} ensures that $u_N$ is convergent in $S'$, and we denote the limit by $u_h=\lim_{N\to\infty} u_N\in S_h^{'}\cap\dot{\mathcal{B}}_{p'}^{-s}$. Recall that  for any $\phi\in S_0\cap\dot{\mathcal{B}}_p^s$, there exists $M\gg 1$ so that $\<u_T,\phi\>=\<u_N,\phi\>$ for any $N>M$. Thus, $\<u_T,\phi\>=\lim_{N\to\infty}\<u_N,\phi\>=\<u_h,\phi\>$, $\forall\phi \in S_0$. Because $S_0\subset S$ is dense in $\dot{\mathcal{B}}_p^s$, we have $u_T=u_h$ as elements of $S'$
			and so $u_T\in \dot{\mathcal{B}}_{p'}^{-s}$. Therefore, the inequality \eqref{eq-eqnorm} shows that the one-to-one mapping $l:				(\dot{\mathcal{B}}_{p}^{s})' 
			\to 
			\dot{\mathcal{B}}_{p'}^{-s}$ is homeomorphic.
			
			Define an equivalent norm $\|\cdot\|_1$ on $\dot{\mathcal{B}}_{p'}^{-s}$ by $\|u\|_1=\sup_{\|\phi\|_{\dot{\mathcal{B}}_p^s}\leq 1} |\<u,\phi\>|$.  Therefore, we have $(\dot{\mathcal{B}}_p^s)'=(\dot{\mathcal{B}}_{p'}^{-s},\|\cdot\|_1)$. Applying the argument above one more time, for all $W\in (\dot{\mathcal{B}}_{p'}^{-s},\|\cdot\|_1)'$, there exists an element $v_W\in \dot{\mathcal{B}}_p^s$ such that $W\psi=\<v_W,\psi\>$, $\forall \psi\in\dot{\mathcal{B}}_{p'}^{-s}$. Thus, $\dot{\mathcal{B}}_p^s$ is reflexive.
			
			(7): By Minkowski inequality and homogeneous Littlewood-Paley theorem, for any $u\in\dot{\mathcal{B}}_p^0$, there exists a function $v\in L_p$ such that $u-v$ is a polynomial. However, $v$ also belongs to $S_h^{'}$ by Remark \ref{rem-realispace}, which tells us that $u=v$. 
		\end{prf}

		\subsection{Bochner spaces}
		Let $E$ be a Banach space. Here we collect some properties of Bochner space $L^p(R^n;E)$, see, e.g., \cite{HNMW}.
		
		For a function $f:\R^n\rightarrow\C$ and an element $e$ of $E$, we denote the tensor product
		$$
		f\otimes e(x)=f(x)e\ .
		$$
		If $V$ is a vector space of complex-valued functions, we use the same notation for the tensor product $$V \otimes E:=\left\{\sum_{n=1}^{N} f_{n} \otimes e_{n}: f_{n} \in V, e_{n} \in E, n=1, \ldots, N ; N=1,2, \ldots\right\}\ .$$
		Now we can introduce some properties of space $L^p(R^n;E)$.
		\begin{prop}[\cite{HNMW}, p.68]\label{prop-density}
			Let $p\in[1,\infty]$,
			the	tensor product $L^p(\R^n)\otimes E$ is dense in $L^p(\R^n;E)$.
		\end{prop}
		\begin{prop}[\cite{HNMW}, p.54]\label{prop-dualthm}
			Let $p\in[1,\infty)$ and $E$ be a reflexive Banach space, then we have
			$$\left(L^{p}(S ; E)\right)'=L^{p'}\left(S ; E'\right)\ ,\  \frac{1}{p}+\frac{1}{p'}=1\ .$$
		\end{prop}


		\end{document}